\documentclass[11pt,hyp,]{amsart}
\usepackage{hyperref}
\hypersetup{nesting=true,debug=true,naturalnames=true}
\usepackage{graphicx,amssymb,upref}

\newtheorem{theorem}{Theorem}
\newtheorem{lemma}[theorem]{Lemma}

\title{Random walk with equidistant multiple function barriers}
\author{Theo van Uem}  
\address{Amsterdam University of Applied Sciences, Amsterdam, The Netherlands.} 
\email{tjvanuem@gmail.com}  

\begin{document}
\hbadness=99999


\begin{abstract}
We obtain expected number of arrivals, absorption probabilities and
expected time before absorption for a discrete random walk on the integers with an
infinite set of equidistant multiple function barriers.
\end{abstract}
\maketitle

\section{Introduction}

Random walk can be used in various disciplines: in medicine and biology where absorbing
barriers give a natural model for a wide variety of phenomena, in physics as a simplified model of
Brownian motion, in ecology to describe individual animal movements and population dynamics.
Random walks have been studied for decades on regular structures such as lattices. Percus \cite{per}
considers an asymmetric random walk, with one or two boundaries, on a one-dimensional lattice.
At the boundaries, the walker is either absorbed or reflected back to the system. Using generating
functions the probability distribution of being at position m after n steps is obtained, as well as the
mean number of steps before absorption. El-Shehawey \cite{el1} \cite{el2} obtains absorption probabilities at
the boundaries for a random walk between one or two partially absorbing boundaries as well as
the conditional mean for the number of steps before stopping given the absorption at a specified
barrier, using conditional probabilities. In this paper we obtain expected number of arrivals,
absorption probabilities and expected time before absorption for a discrete random walk on the
integers with an infinite set of equidistant multiple function barriers. A multiple function barrier
(MFB) is a state that can absorb, reflect, let through or hold for a moment. In each mfb we have
probabilities $p_0,q_0,r_0,s_0$ for moving forward and backward , staying for a moment in the MFB
and absorption in the MFB, where $p_0+q_0+r_0+s_0=1, \ p_0q_0s_0>0$. MFB’s of type $p_0q_0r_0s_0$ are
defined in each barrier $kN(\ k\in \mathbb Z,\ N>1)$.
The random walk between the MFB’s is of $pqr$ type, where $p$ is the one-step forward probability, $q$
one-step backward probability $(pq>0)$ and $r=1-p-q$ the probability to stay for a moment in the same position. 
We start in $i_0 \ (0\leq i_0<N)$.
\section{Random walk on the MFB's}
We define the expected number of arrivals in state $j$ when starting in state $i$:
\[
x_j=x_{i,j}=\sum_{k=0}^{\infty}p_{ij}^{(k)}
\]
Let $
\rho=\frac{p}{q} $ and $\lambda_1$ and $\lambda_2 (\lambda_1\geq \lambda_2)$ are the solutions of $q\lambda^2-(1-r)\lambda+p=0$. If $p>q$ then $\lambda_1=\rho,\lambda_2=1$. If $p<q$ then $\lambda_1=1,\lambda_2=\rho$. If $p=q$ then $\lambda_1=\lambda_2=1$.
We start  with a $pqr$ random walk on the integers:
\begin{lemma} \label{t1}
\begin{equation} \label{1}
x_{n}=\delta(n,i_0)+px_{n-1}+qx_{n+1}+rx_{n} \quad (n\in 
\mathbb Z)\quad \rho\neq 1
\end{equation}
 has solution:
\begin{equation}
x_{n}=\left\{\begin{array}{l}\frac{\lambda_1^{n-i_0}}{\sqrt{(1-r)^2-4pq}}\ \ \ \ \ \ \ \    (n\leq i_0) \\
\frac{\lambda_2^{n-i_0}}{\sqrt{(1-r)^2-4pq}}\ \ \ \ \ \ \ \    (n\geq i_0) \\\end{array}\right.
\end{equation}
\end{lemma}
\begin{proof}
Let 
\begin{equation*}
G(s)=\sum_{k=-\infty}^{\infty}x_{k}s^k \quad(|s|<1)
\end{equation*}
Using \ref{1} we obtain:
\[G(s)=s^{i_0}+psG(s)+qs^{-1}G(s)+rG(s)\]
\[G(s)=\frac{s^{i_0}}{1-ps-qs^{-1}-r}\]
We use the inverse z-transform: $x_{n}=\frac{1}{2\pi \mathrm i}\oint H(z)z^{n-1} \mathrm d z$, where the integration is along the circle $|z|=1$ and anticlockwise. we have:
\begin{equation*}
H(z)=\sum_{n=-\infty}^{\infty}f_{n}z^{-n}=G(z^{-1})=
\frac{z^{-i_0}}{1-pz^{-1}-qz-r}.
\end{equation*}
So,
\[x_{n}=\frac{1}{2\pi \mathrm i}\oint \frac{z^{n-1-{i_0}}}{1-pz^{-1}-qz-r}=\frac{1}{2\pi \mathrm i}\oint \frac{-z^{n-{i_0}}}{q(z-\lambda_1)(z-\lambda_2)}  \mathrm d z
\]
Apply the residue theorem.
\end{proof}

\begin{theorem}
The random walk on the subset of equidistant MFB's is described by the difference equations:
CASE $\rho\neq 1 $.
\begin{equation*}
(\lambda_1-\lambda_2)q_0x_{(k+1)N}+\omega_0x_{kN}+(\lambda_1-\lambda_2)p_0\rho^{N-1}x_{(k-1)N}=
\end{equation*}
\begin{equation}\label{3}
(\lambda_2^{N-i_0}-\lambda_1^{N-i_0})\delta(k,0)+(\lambda_1^{-i_0}-\lambda_2^{-i_0})\rho^N\delta(k,1)\quad (k\in \mathbb Z)
\end{equation}
where
\begin{equation} \label{4}
\omega_0=(\lambda_2^{N}-\lambda_1^{N})(1-r_0)+(\lambda_1^{N-1}-\lambda_2^{N-1})(p_0+q_0\rho)
\end{equation}
CASE $\rho= 1$.
\begin{equation}
q_0x_{(k+1)N}-(p_0+q_0+Ns_0)x_{kN}+p_0x_{(k-1)N}=
(i_0-N)\delta(k,0)-i_0\delta(k,1)\quad 
\end{equation}
\end{theorem}

\begin{proof}
We start with $0<i_0<N$. Random walk on interval $[kN+1,(k+1)N-1]$:
\begin{equation}\label{5}
(1-r)x_{kN+n}=\delta(k,0)\delta(n,i_0)+px_{kN+n-1}+qx_{kN+n+1} \quad (n=2,3,\dots,N-2)
\end{equation}
Characteristic equation:
\[q\lambda^2-(1-r)\lambda+p=0\]
A general solution of \ref{5} is (use Lemma \ref{t1}):
\begin{equation} \label{6}
x_{kN+n}=\left\{\begin{array}{l}\frac{\lambda_1^{n-i_0}\delta(k,0)}{\sqrt{(1-r)^2-4pq}}+a_k\lambda_1^n+b_k\lambda_2^n\ \ \ \ \ \ \ \    (n=1,\dots,i_0) \\
\frac{\lambda_2^{n-i_0}\delta(k,0)}{\sqrt{(1-r)^2-4pq}}+a_k\lambda_1^n+b_k\lambda_2^n\ \ \ \ \ \ \ \    (n=i_0,\dots,N-1) \\\end{array}\right.
\end{equation}
Let $\zeta=[(1-r)^2-4pq]^{-\frac{1}{2}}$. By focusing on states $kN+1$ and $(k+1)N-1$ we get:
\[x_{kN+1}=p_0x_{kN}+qX_{kN+2}+rx_{kN+1} \]
\[x_{(k+1)N-1}=px_{(k+1)N-2}+q_0x_{(k+1)N}+rx_{(k+1)N-1} \]

\[p_0x_{kN}=p[\zeta\lambda_1^{-i_0}\delta(k,0)+a_k+b_k]\]
\[ q_0x_{(k+1)N}=q[\zeta\lambda_2^{N-i_0}\delta(k,0)+a_k\lambda_1^N+b_k\lambda_2^N]\]

\[(\lambda_2^N-\lambda_1^N)a_k=\lambda_2^N\frac{p_0}{p}x_0-\frac{q_0}{q}x_N+\zeta\lambda_2^N(\lambda_2^{-i_0}-\lambda_1^{-i_0})\delta(k,0)\]
\[(\lambda_2^N-\lambda_1^N)b_k=-\lambda_1^N\frac{p_0}{p}x_0+\frac{q_0}{q}x_N+\zeta\lambda_1^{N-i_0}-\lambda_2^{N-i_0})\delta(k,0)\]

Focusing on state $kN$:
\[x_{kN}=px_{kN-1}+qx_{kN+1}+r_0x_{kN}\]
After some calculations, we get \ref{3}

CASE $\rho=1$
We use the same method, where (verified by substitution):
\begin{equation} \label{7}
x_{kN+n}=\left\{\begin{array}{l}a_kn+b_k+\frac{n-i_0}{p}\ \ \ \ \ \ \ \    (n=1,\dots,i_0) \\
a_kn+b_kn\ \ \ \ \ \ \ \    (n=i_0,\dots,N-1) \\\end{array}\right.
\end{equation}
The special case where we start in $i_0=0$ can be handled in the same way, resulting in \ref{3} and \ref{5} with $i_0=0$ when $\rho\neq 1$ respectively $\rho=1$.
\end{proof}
\begin{theorem}
The RW on the MFB's is symmetric if and only if $(i_0=0) \land (q_0=p_0\rho^{N-1})$
\end{theorem}
\begin{proof}
See \ref{3} and \ref{5}.
\end{proof}
Notice that  $p_0p^{N-1}=q_0q^{N-1}$ can be interpreted as: direct probability from a MFB to it's right neighbor equals direct probability in the reverse direction.
\section{Value of the MFB game}
We define a moment generating function on the MFB's:
\begin{equation}
F(s)=\sum_{k=-\infty}^{\infty}x_{kN}s^k \quad(|s|<1)
\end{equation}
\begin{theorem}
CASE $\rho\neq 1$:
\begin{equation}
F(s)=\frac{\lambda_2^{N-i_0}-\lambda_1^{N-i_0}+(\lambda_1^{-i_0}-\lambda_2^{-i_0})\rho^N s}{(\lambda_1-\lambda_2)q_0s^{-1}+\omega_0+(\lambda_1-\lambda_2)p_0\rho^{N-1}s}
\end{equation}
CASE $\rho= 1$:
\begin{equation}
F(s)=\frac{i_0-N-i_0s}{q_0s^{-1}-(p_0+q_0+Ns_0)+p_0s}
\end{equation}
\end{theorem}
\begin{proof}
Use \ref{3} and \ref{5}.
\end{proof}
\begin{theorem}
Probability of absorption in a MFB is 1:
\[\sum_{k=-\infty}^{\infty}s_0x_{kN}=1\]
\end{theorem}
\begin{proof}
In both cases we have: $F(1)=\sum_{k=-\infty}^{\infty}x_{kN}=\frac{1}{s_0}$
\end{proof}
We define the value $v$ of the MFB game as: $v=\sum_{k=-\infty}^{\infty}kx_{kN}$.
\begin{theorem}
CASE $\rho\neq 1$:
\begin{multline*}
v=\frac{(\lambda_1^{-i_0}-\lambda_2^{-i_0})\rho^N}{(\lambda_2^N-\lambda_1^N)s_0}+ \\
-\frac{(\lambda_1-\lambda_2)(q_0-p_0\rho^{N-1})[(\lambda_1^{-i_0}-\lambda_2^{-i_0})\rho^N+\lambda_2^{N-i_0}-\lambda_1^{N-i_0}]}{(\lambda_2^N-\lambda_1^N)^2s_0^2}
\end{multline*}
CASE $\rho=1$:
\[
v=\frac{p_0-q_0+i_0s_0}{Ns_0^2}
\]
\end{theorem}
\begin{proof}
$v=[\frac{\mathrm{d}F}{\mathrm{d}s}]_{s=1}.$
\end{proof}
Notice that the symmetric random walk on the MFB's has value 0.
\section{Expected number of arrivals}
\begin{theorem} \label{t4}
The expected number of arrivals to the MFB's is:
CASE $\rho\neq 1$:
\begin{equation}
x_{kN}=\left\{\begin{array}{l}\{(\lambda_1^{N-i_0}-\lambda_2^{N-i_0})\xi_1+\rho^N(\lambda_2^{-i_0}-\lambda_1^{-i_0})\}\Omega \xi_1^{k-1}\ \ \ \ \ \ \ \    (k\leq 0) \\
\{(\lambda_1^{N-i_0}-\lambda_2^{N-i_0})\xi_2+\rho^N(\lambda_2^{-i_0}-\lambda_1^{-i_0})\}\Omega \xi_2^{k-1}\ \ \ \ \ \ \ \    (k\geq 1)
 \\
 \end{array}\right.
\end{equation}
where
\begin{equation}
(\lambda_1-\lambda_2)q_0\xi_i^2+\omega_0\xi_i+(\lambda_1-\lambda_2)p_0\rho^{N-1}=0 \quad (i=1,2) \quad \xi_1>1>\xi_2>0
\end{equation}
\[\Omega=\{\omega_0^2-4p_0q_0(1-\rho)^2\rho^{N-1}\}^{-\frac{1}{2}}\]

CASE $\rho=1$:
\begin{equation}
x_{kN}=\left\{\begin{array}{l}\frac{\{(N-i_0)\xi_1+i_0\}\xi_1^{k-1}}{\sqrt{(p_0+q_0+Ns_0)^2-4p_0q_0}}\ \ \ \ \ \ \ \    (k\leq 0) \\
\frac{\{(N-i_0)\xi_2+i_0\}\xi_2^{k-1}}{\sqrt{(p_0+q_0+Ns_0)^2-4p_0q_0}}\ \ \ \ \ \ \ \    (k\geq 1)
 \\
 \end{array}\right.
\end{equation}
where
\begin{equation}
q_0\xi_i^2-(p_0+q_0+Ns_0)\xi_i+p_0=0 \quad (i=1,2) \quad \xi_1>1>\xi_2>0
\end{equation}
\end{theorem}
\begin{proof}
CASE $\rho\neq 1$
We use the inverse z-transform: \\ $x_{kN}=\frac{1}{2\pi \mathrm i}\oint H(z)z^{k-1} \mathrm d z$, where the integration is along the circle $|z|=1$ and anticlockwise. Using \ref{3} we get:
\begin{equation*}
H(z)=\sum_{n=-\infty}^{\infty}f_{n}z^{-n}=F(z^{-1})=
\frac{(\lambda_2^{N-i_0}-\lambda_1^{N-i_0})z+(\lambda_1^{-i_0}-\lambda_2^{-i_0})\rho^n }{(\lambda_1-\lambda_2)q_0z^2+\omega_0z+(\lambda_1-\lambda_2)p_0\rho^{N-1}}
\end{equation*}
So,
\[x_{n}=\frac{1}{2\pi \mathrm i}\oint \frac{(\lambda_2^{N-i_0}-\lambda_1^{N-i_0})z^n+(\lambda_1^{-i_0}-\lambda_2^{-i_0})\rho^n z^{n-1} }{(\lambda_1-\lambda_2)q_0(z-\xi_1)(z-\xi_2)}  \mathrm d z
\]
Apply the residue theorem.
CASE $\rho= 1$
Using \ref{5} we get:
 \begin{equation}
F(z^{-1})=\frac{i_0-N-i_0z^{-1}}{q_0z-(p_0+q_0+Ns_0)+p_0z^{-1}}=\frac{(i_0-N)z-i_0}{q_0(z-\xi_1)(z-\xi_2)}
\end{equation}
Use $x_{kN}=\frac{1}{2\pi \mathrm i}\oint F(z^{-1})z^{k-1} \mathrm d z$ and the residue theorem.
\end{proof}
\begin{theorem}
CASE $\rho\neq 1$:
\[(1-\rho^N)x_{kN+n}= \]
\begin{equation} 
\left\{\begin{array}{l}\frac{p_0}{p}[\rho^{n-kN}-\rho^N]x_{kN}+\frac{q_0}{q}[1-\rho^{n-kN}]x_{(k+1)N} +\frac{(1-\rho^n)(\rho^{N-i_0}-1)}{p-q}\delta(k,0) \\
(n=1,\dots,i_0) \\
\frac{p_0}{p}[\rho^{n-kN}-\rho^N]x_{kN}+\frac{q_0}{q}[1-\rho^{n-kN}]x_{(k+1)N} +\frac{(\rho^n-\rho^N)(1-\rho^{-i_0})}{p-q}\delta(k,0)\\
(n=i_0,\dots,N-1) \\\end{array}\right.
\end{equation}
CASE $\rho=1$:
\begin{equation}
x_{kN+n}=
\left\{\begin{array}{l}\frac{p_0(N-n)x_{kN}+q_0nx_{(k+1)N}+n(N-i_0)\delta(k,0)}{pN}\ \ \ \ \ \ \ \    (n=1,\dots,i_0) \\
\frac{p_0(N-n)x_{kN}+q_0nx_{(k+1)N}+i_0(N-n)\delta(k,0)}{pN} \ \ \ \ \ \ \ \    (n=i_0,\dots,N-1) \\\end{array}\right.
\end{equation}
\end{theorem}
\begin{proof}
Along the same lines as in Theorem \ref{t4}, using \ref{6} and \ref{7}.
\end{proof}

\section{Mean absorption time}
Let $m_i$ be the mean absorption time (in any MFB) when starting in state $i$ where $i \in \mathbb Z$.
\begin{theorem} \label{t5}
\[ m_i=m_{i\bmod N} \quad (i \in \mathbb Z)\]
CASE $\rho\neq 1$. If $0\leq i\leq N$:
\begin{multline*}
m_i=\frac{N\rho^{-i}}{(q-p)(1-\rho^{-N})}+\frac{i}{q-p}+\frac{1}{s_0}+ \\
\frac{p_0+q_0(N-1)}{(q-p)s_0}+\frac{N[p_0\rho^{-1}+q_0\rho^{1-N}+r_0-1]}{(q-p)(1-\rho^{-N})s_0}
\end{multline*}
CASE $\rho= 1$. If $0\leq i\leq N$:
\[
m_i=\frac{i(N-i)}{2p}+\frac{1}{s_0}+\frac{p_0+q_0(N-1)}{2ps_0}
\]
\end{theorem}
\begin{proof}
\[m_i=p(m_{i+1}+1)+q(m_{i-1}+1)+r(m_i+1)\quad (1\leq i\leq N-1)\]
\[m_0=p_0(m_{1}+1)+q_0(m_{-1}+1)+r_0(m_0+1)+s_0.1\]
Because of \[ m_i=m_{i\bmod N} \quad (i \in \mathbb Z)\] we have:
\begin{equation}\label{10}
(1-r)m_i=pm_{i+1}+qm_{i-1}+1\quad (1\leq i\leq N-1)
\end{equation}
\[m_0=m_N\]
\[(1-r_0)m_0=p_0m_{1}+q_0(m_{N-1}+1)+1\]
Use $m_i=a\rho^{-i}+b+\frac{i}{q-p}$ (case $\rho\neq 1$) or $m_i=ai+b-\frac{i^2}{2p}$ (case $\rho= 1$) where $0\leq i\leq N$ because $m_0$ and $m_N$ are part of the difference pattern \ref{10}.
\end{proof}
Notice that the results for $\rho=1$ can also be obtained by applying l'Hospitals rule in the result for $\rho\neq 1$ (except Theorem \ref{t5} where we need l'Hospitals rule twice).

\end{document}